\documentclass{amsart}

\usepackage{amssymb, amsmath}
\usepackage{mathrsfs}
\usepackage{amscd}
\usepackage[active]{srcltx}
\usepackage{verbatim}

\usepackage[colorlinks,linkcolor={blue},
citecolor={blue},urlcolor={red},]{hyperref}

\usepackage{newsymbol}
\let\mathcal \undefined
\def\mathcal{\mathscr}

\let\emptyset \undefined
\let\ge       \undefined
\let\le       \undefined
\newsymbol\le          1336  
\newsymbol\ge          133E  
\newsymbol\emptyset    203F
\newsymbol\notle       230A
\newsymbol\notge       230B

\theoremstyle{plain}
\newtheorem{theorem}{Theorem}[section]
\theoremstyle{remark}
\newtheorem{remark}[theorem]{Remark}
\newtheorem{example}[theorem]{Example}
\theoremstyle{plain}
\newtheorem{corollary}[theorem]{Corollary}
\newtheorem{lemma}[theorem]{Lemma}
\newtheorem{proposition}[theorem]{Proposition}
\newtheorem{definition}[theorem]{Definition}

\numberwithin{equation}{section}

\def\Z{{\mathbb Z}}

\def\R{{\mathbb R}}
\def\C{{\mathbb C}}
\def\K{{\mathbb K}}

\renewcommand{\P}{{\mathbb P}}

\newcommand{\BA}{\mathscr{BA}}

\renewcommand{\a}{\alpha}

\newcommand{\e}{\varepsilon}

\newcommand{\om}{\omega}
\newcommand{\Om}{\Omega}

\renewcommand{\L}{L^2(0,T)}

\def\L{\calL}

\def\M{{\mathscr M}}
\def\parti{{\mathscr P}}

\def\T{\mathbb{T}}

\def\wot{_{\rm w}}
\def\sot{_{\rm s}}

\def\ba{\begin{eqnarray*}}
\def\ea{\end{eqnarray*}}

\def\bbe{\begin{equation}}

\def\ee{\end{equation}}

\newcommand{\beq}{\begin{equation}}
\newcommand{\eeq}{\end{equation}}
\newcommand{\bal}{\begin{aligned}}
\newcommand{\eal}{\end{aligned}}
\newcommand{\ben}{\begin{enumerate}}
\newcommand{\een}{\end{enumerate}}
\newcommand{\bit}{\begin{itemize}}
\newcommand{\eit}{\end{itemize}}

\newcommand{\bth}{\begin{theorem}}
\renewcommand{\eth}{\end{theorem}}
\newcommand{\bpr}{\begin{proposition}}
\newcommand{\epr}{\end{proposition}}
\newcommand{\ble}{\begin{lemma}}
\newcommand{\ele}{\end{lemma}}
\newcommand{\bpf}{\begin{proof}}
\newcommand{\epf}{\end{proof}}
\newcommand{\bex}{\begin{example}}
\newcommand{\eex}{\end{example}}
\newcommand{\bre}{\begin{example}}
\newcommand{\ere}{\end{example}}

\newcommand{\calL}{{\mathcal L}}
\newcommand{\n}{\Vert}

\newcommand{\embed}{\hookrightarrow}
\newcommand{\s}{^*}
\newcommand{\lb}{\langle}
\newcommand{\rb}{\rangle}

\newcommand{\limn}{\lim_{n\to\infty}}

\begin{document}

\title[Spaces of operator-valued functions]
{Spaces of operator-valued functions measurable with respect to the strong operator topology}
\author{Oscar Blasco}
\address{University of Valencia\\Departamento de An\'alisis Matem\'atico\\46100 Burjassot, Valencia\\Spain}
\email{oscar.blasco@uv.es}
\author{Jan van Neerven}
\address{Delft University of Technology\\Delft Institute of Applied Mathematics\\P.O. Box 5031\\2600 GA Delft\\The Netherlands}
\email{J.M.A.M.vanNeerven@TUDeflt.nl}
\thanks{The first named author is partially supported by the spanish project MTM2008-04594/MTM. The second named author is supported by VICI subsidy 639.033.604 of the Netherlands Organisation for Scientific Research (NWO)}
\keywords{Operator-valued functions, operator-valued multipliers,
vector meas\-ures} \subjclass[2000]{28B05, 46G10}

\begin{abstract}
Let $X$ and $Y$ be Banach spaces and
$(\Om, \Sigma,\mu)$ a finite measure space. In this note we introduce the space $L^p[\mu; \L(X,Y)]$ consisting of all (equivalence classes of) functions $\Phi:\Om \mapsto \L(X,Y)$ such that $\om\mapsto \Phi(\om)x$ is strongly $\mu$-measurable for all $x\in X$ and
$\om\mapsto \Phi(\om)f(\om)$ belongs to $L^1(\mu;Y)$ for all $f\in L^{p'}(\mu;X)$, $1/p+1/p'=1$. We show that functions in $L^p[\mu;\L(X,Y)]$ define operator-valued measures with bounded $p$-variation and use
these spaces to obtain an isometric characterization of the space of all $\calL(X,Y)$-valued
multipliers acting boundedly from $L^p(\mu;X)$ into $L^q(\mu;Y)$,
$1\le q< p<\infty$.
\end{abstract}

\maketitle

\section{Introduction}
Let $(\Om,\Sigma,\mu)$ be a finite measure space and let $X$ and $Y$ be Banach spaces
over $\K=\R$ or $\C$.
In his talk at the 3rd meeting on Vector Measures, Integration and Applications (Eichst\"att, 2008), Jan Fourie presented some applications
of the following extension of an elementary observation due to Bu and Lin \cite[Lemma 1.1]{BuLin}.

\begin{proposition}\label{fourie}
Let $\Phi:\Om\to \calL(X,Y)$ be a strongly $\mu$-measurable function. For all $\e>0$ there exists strongly $\mu$-measurable function
$f^\e:\Om\to X$ such that for $\mu$-almost all $\om\in \Om$ one has $\n f^\e(\om)\n \le 1$ and
$$ \n \Phi(\om)\n \le \n \Phi(\om)f^\e(\om)\n + \e.$$
\end{proposition}

Recall that a function $\phi:\Om\to Z$, where $Z$ is a Banach space, is said to be {\em strongly $\mu$-measurable} if there exists a sequence of $\Sigma$-measurable simple functions $\phi_n:\Om\to Z$ such that for $\mu$-almost all $\om\in\Om$ one has $\limn \phi_n(\om) = \phi(\om)$ in $Z$.

In Proposition \ref{fourie}, the strong $\mu$-measurability
assumption on $\Phi$ refers to the norm of $\calL(X,Y)$ as a
Banach space. The next two examples show that the conclusion
of Proposition \ref{fourie} often holds if we impose merely strong $\mu$-measurability of the orbits of $\Phi$.

\begin{example} \label{ex1} Consider $X=\ell^\infty(\Z)$, let $\T$ be the unit circle, and define $\Phi:\T\to\ell^\infty(\Z) = \calL(\ell^1(\Z),\K)$ by $\Phi(t):=(e^{int})_{n\in \Z}$. For all $x\in \ell^1(\Z)$
the function $t\mapsto \Phi(t)x = \sum_{n\in\Z} x_n e^{int}$ is
continuous, but the function  $t\mapsto \Phi(t)$ fails to be strongly measurable. Taking for $f$ the constant function with value $u_0\in
\ell^1(\Z)$, defined by $u_0(0) =1$ and $u_0(n)=0$ for $n\not=0$,
we have
$$\|\Phi(t)\|= |\Phi(t)f(t)| =  |\langle u_0,\Phi(t)\rangle|=1\quad \forall t\in \T.$$
 \end{example}

 \begin{example} \label{ex2} \ Consider \, $X=C([0,1])$ \, and \, define \, $\Phi:[0,1]\to M([0,1])= \calL(C([0,1]),\K)$\,
by $\Phi(t):=\delta_t$. For all $x\in X$ the function
$t\mapsto \Phi(t)x = x(t)$ is continuous, but the function  $t\mapsto \Phi(t)$
fails to be strongly measurable. If $f:[0,1]\to X$ is a strongly measurable function such that $(f(t))(t)=1$ for all $t\in [0,1]$ (e.g., take $f(t)\equiv 1)$, we have
$$\|\Phi(t)\|= |\langle f(t),\Phi(t)\rangle|=1\quad \forall t\in [0,1].$$
 \end{example}

Thus it is natural to ask whether strong $\mu$-measurability of $\Phi$ can
be weakened to strong $\mu$-measurability of the orbits $\om\mapsto \Phi(\om)x$ for all $x\in X$, or even to $\mu$-measurability of the functions $\om\mapsto \|\Phi(\om)x\|$. Although in general the
answer is negative even when $\dim Y=1$ (Example \ref{ex}),
various positive results can be formulated under additional
assumptions on $X$ or $\Phi$ (Propositions \ref{prop}, \ref{prop:BA}, and their corollaries).

One of the applications of Proposition \ref{fourie} was the study of multipliers between
spaces of vector-valued integrable functions. In  \cite{FouSch}, for $1\le p,q<\infty$, $\hbox{Mult\/}(L^p(\mu;X),L^q(\mu;Y))$ is defined to be the space of all strongly $\mu$-measurable functions $\Phi:\Omega \mapsto\L(X,Y)$ such that $\om\mapsto\Phi(\om)f(\om)$ belongs to
$L^q(\mu;Y)$ for all $f\in L^p(\mu;X)$. It is shown (see \cite[Proposition 3.4]{FouSch}) that for $1\le q<p<\infty$ and $1/r=1/q-1/p$ one has a natural isometric isomorphism $$\hbox{Mult\/}(L^p(\mu;X),L^q(\mu;Y))\simeq  L^r(\mu;\L(X,Y)).$$
We observe (Proposition \ref{meas}) that the strong $\mu$-measurability of $\Phi$ as function with values in $\L(X,Y)$ is not really needed to define bounded operators from $L^p(\mu;X)$ into $L^q(\mu;Y)$; it is possible
to weaken the measurability assumptions on the multiplier functions by only requiring strong $\mu$-measurability of its orbits. This will motivate the introduction of an intermediate space between $L^p(\mu; \L(X,Y))$ and the space $L\sot^p(\mu; \L(X,Y))$ of functions $\Phi:\Omega\mapsto \L(X,Y)$ such that $\om \to \Phi(\om)x$ belongs to $L^p(\mu;Y)$ for all $x\in X$.  This is done by selecting the functions in $L\sot^p(\mu; \L(X,Y))$ for which $\om\mapsto\Phi(\om)f(\om)$ belongs to $L^1(\mu;Y)$ for all $f\in L^{p'}(\mu;X)$, $1/p+1/p'=1$. We shall denote this space by $L^p[\mu; \L(X,Y)]$.  We shall see that,  for $1\le p<\infty$, functions in this space define $\L(X,Y)$-valued measures of bounded $p$-variation (Theorems \ref{t2} and \ref{t1}),
and prove that one has a natural isometric isomorphism $$\hbox{Mult\/}[L^p(\mu;X),L^q(\mu;Y)]\simeq  L^r[\mu;\L(X,Y)],$$
where $1/r=1/q-1/p$ and $\hbox{Mult\/}[L^p(\mu;X),L^q(\mu;Y)]$ is defined to be the linear space of all functions $\Phi:\Omega \mapsto\L(X,Y)$ such that  $\om\mapsto\Phi(\om)x$ is strongly $\mu$-measurable for all $x\in X$ and $\om\mapsto\Phi(\om)f(\om)$ belongs to
$L^q(\mu;Y)$ for all $f\in L^p(\mu;X)$ (Theorem \ref{thm:mult}).

\section{Strong $\mu$-normability of operator-valued functions}

Let $(\Om,\Sigma,\mu)$ be a finite measure space and let $X$ and $Y$ be Banach spaces. 

\begin{definition} Consider a function $\Phi:\Om\to \calL(X,Y)$.
\ben
\item $\Phi$ is called {\em strongly
$\mu$-normable} if for all $\e>0$ there exists strongly
$\mu$-measurable function $f^\e:\Om\to X$ such that for
$\mu$-almost all $\om\in \Om$ one has $\n f^\e(\om)\n \le 1$ and
$$ \n \Phi(\om)\n \le \n \Phi(\om)f^\e(\om)\n + \e.$$
\item $\Phi$ is called {\em weakly
$\mu$-normable} if for all $\e>0$ there exist strongly
$\mu$-measurable functions $f^\e:\Om\to X$ and $g^\e:\Om\to Y\s$ such that for
$\mu$-almost all $\om\in \Om$ one has $\n f^\e(\om)\n \le 1$, $\n g^\e(\om)\n\le 1$, and
$$ \n \Phi(\om)\n \le |\lb \Phi(\om)f^\e(\om),g^\e(\om)\rb| + \e.$$
\een
\end{definition}

Clearly, every weakly $\mu$-normable function is strongly $\mu$-normable.
In the case $Y=\K$ the notions of weak and strong $\mu$-normability coincide and we shall simply speak of {\em normable} functions. 

It will be convenient to formulate our results on $\mu$-normability in the following more general setting.
Let $S$ an arbitrary nonempty set. A function $f:\Om\to S$ is called a {\em $\Sigma$-measurable elementary function}
if for $n\ge 1$ there exist disjoint sets $A_n\in \Sigma$ and elements $s_n\in S$ such that $\bigcup_{n\ge 1} A_n = \Om$ and
$f = \sum_{n\ge 1} 1_{A_n}\otimes s_n$.
Since no addition is defined in $S$, this sum should be interpreted as shorthand notation
to express that $f\equiv s_n$ on $A_n$.
A function $g:S\to \R$ is called {\em bounded from above} if
$\sup_{s\in S} \,g(s) <\infty.$ The set of all such functions is denoted by $\BA(S)$.

\begin{proposition}\label{prop}
Let $\Phi: \Om\to \BA(S)$ be such that for all $s\in S$ the
function $ \om\mapsto(\Phi(\om))(s)$ is $\mu$-measurable. If there
is a countable subset $C$ of $S$ such that for all $\phi\in
\Phi(\Omega)$ we have
$$ \sup_{s\in S} \phi(s) = \sup_{s\in C} \phi(s),$$
then for all $\e>0$ there exists a $\Sigma$-measurable elementary function $f^\e:\Om\to S$ such that for $\mu$-almost all $\om\in\Om$ one has
$$\sup_{s\in S}\, (\Phi(\om))(s) \le (\Phi(\om))(f^\e(\om))+\e.$$
\end{proposition}
\begin{proof} The function
$\om\mapsto \sup_{s\in C}(\Phi(\om))(s)$ is $\mu$-measurable, as it is the pointwise supremum of
a countable family of $\mu$-measurable functions.
Let $(s^{(n)})_{n\ge 1}$ be an enumeration of $C$.
For $n\ge 1$ put
$$A_{n} := \big\{\om\in\Om: \ \sup_{s\in S} \Phi(\om)(s) \le  (\Phi(\om))(s^{(n)}) +\e\big\}.$$
These sets are $\mu$-measurable, and therefore there exist sets $A_n'\in\Sigma$ such that $\mu(A_{n}\Delta A_{n}')=0$. Also, $\bigcup_{n\ge 1}A_n = \Om$.
Put $B_1:= A_{1}'$ and $B_{n+1}:= A_{n+1}'\setminus \bigcup_{m=1}^n B_n$ for $n\ge 1$. The sets $B_n$ are $\Sigma$-measurable, disjoint.
Since $B_0:= \Om\setminus\bigcup_{n\ge 1} B_n$ is a $\mu$-null set in $\Sigma$, the function
$$ f^\e := \sum_{n\ge 0} 1_{B_n}\otimes s^{(n)},$$
where $s^{(0)}\in S$ is chosen arbitrarily, has the desired properties.
\end{proof}

From this general point of view one obtains the following corollary.

\begin{corollary}\label{cor1}
Let $X$ and $Y$ be Banach spaces and consider a function $\Phi: \Om\to \calL(X,Y)$.

\ben
\item If $X$ is separable and $\om\mapsto \n \Phi(\om)x\n$ is $\mu$-measurable for all $x\in X$, then $\Phi$ is strongly $\mu$-normable;
\item If $X$ and $Y$ are separable and $\om\mapsto |\lb \Phi(\om)x,y\s\rb|$ is $\mu$-measurable for all $x\in X$ and $y\s\in Y\s$, then $\Phi$ is weakly $\mu$-normable.
\een
\end{corollary}
\begin{proof}
To prove (2) we apply Proposition \ref{prop} to the set $S = B_{X\times Y\s}$ (the unit ball of $X\times Y\s$ with respect to the norm $\n (x,y\s)\n = \max\{\n x\n,\n y\s\n\}$) and the functions $\om\mapsto |\lb \Phi(\om)x, y\s\rb|$, and note that $\Sigma$-measurable elementary functions with values in a Banach space are strongly $\mu$-measurable. Since $X$ is separable, for $C$ we may take a set of the form $\{(x_j,y_k\s): j,k\ge 1\}$, where $(x_j)_{j\ge 1}$
is a dense sequence in $B_X$ and $(y_k\s)_{k\ge 1}$ is a sequence in $B_{Y\s}$ which is norming for $Y$. 

The proof of (1) is similar.
\end{proof}

\begin{proof}[Proof of Proposition \ref{fourie}]
By assumption, $\Phi$ can be approximated $\mu$-almost everywhere by a sequence of simple functions with values in $\L(X,Y)$. Each one of the countably many operators in the ranges of these functions is normed by some separable subspace of $X$.
This produces a separable closed subspace $\widetilde X$ of $X$ such that
for $\mu$-almost all $\om\in \Om$,
$$\n \Phi(\om)\n_{\calL(X,Y)} = \n \Phi(\om)\n_{\calL(\widetilde X,Y)}.$$
Now we may apply Corollary \ref{cor1}(1).
\end{proof}

Instead of a countability assumption on the set $S$ we may also impose regularity assumptions on $\mu$ and $\Phi$:

\begin{proposition}\label{prop:BA}
Let $\mu$ be a finite Radon measure on a topological space $\Om$.
Let $\Phi: \Om\to \BA(S)$ be such that for all $s\in S$
the function $ \om\mapsto(\Phi(\om))(s)$ is lower semicontinuous.
Then for all $\e>0$ there exists a Borel measurable elementary function $f^\e:\Om\to S$ such that for $\mu$-almost all $\om\in\Om$ one has
$$\sup_{s\in S}\, (\Phi(\om))(s) \le (\Phi(\om))(f^\e(\om))+\e.$$
\end{proposition}
\begin{proof}
Let us first note that the function $$m(\om) := \sup_{s\in S} (\Phi(\om))(s)$$ is lower semicontinuous, since it is the pointwise supremum of a family of lower semicontinuous functions. In particular, $m$ is Borel measurable.

Fix $\e>0$.
Using Zorn's lemma, let $(\Om_i)_{i\in I}$ be a maximal collection of disjoint Borel sets such that the following two properties are satisfied for all $i\in I$:
\ben
\item[(a)] $\mu(\Om_i)>0$;
\item[(b)] there exists $s_i\in S$ such that $m(\om)\le (\Phi(\om))(s_i) +\e$ for all $\om\in \Om_i$.
\een
Clearly, (a) implies that the index set $I$ is countable. We claim that $$\mu\Big(\Om\setminus\bigcup_{i\in I} \Om_i\Big) =0.$$
The proof is then finished by taking $f^\e := \sum_{i\in I} 1_{\Om_i}\otimes s_i$ and extending this definition to the remaining Borel $\mu$-null set
by assigning an arbitrary constant value on it; by (b) and the claim, this function satisfies the required inequality $\mu$-almost everywhere.

To prove the claim let $\Om' := \Om\setminus\bigcup_{i\in I} \Om_i$
and suppose, for a contradiction, that $\mu(\Om')>0$.
By passing to a Borel subset of $\Om'$ we may assume that
$\sup_{\om'\in\Om'}m(\om')<\infty$.
Let $$M := \hbox{ess sup}_{\om'\in \Om'} \ m(\om').$$
The set $$A := \{\om'\in\Om': \ m(\om') \ge M - \tfrac13\e\}$$
is Borel and
satisfies $\mu(A)>0$.
Since $\mu$ is a Radon measure we may select a compact set $K$ in $\Om$ such that $K\subseteq A$ and $\mu(K)>0$.
For any $\om'\in K$ we can find $s'\in S$ such that
$$ m(\om') \le  (\Phi(\om'))(s') +\tfrac13\e. $$ By lower semicontinuity, the set
$$O' :=\big\{\om\in \Om: \ (\Phi(\om'))(s') < (\Phi(\om))(s') +\tfrac13\e\big\}$$
is open and contains $\om'$. Choosing such an open set for every $\om'\in K$, we obtain an open cover of $K$, which therefore has a finite subcover. At least one of the finitely many open sets of this subcover intersects $K$ in a set of positive measure. Hence, there exist $\om_0\in K$ and $s_0\in S$, as well as an
open set $O_0\subseteq \Om$ such that $\om_0\in O_0$,  $\mu(K\cap O_0)>0$,
$$ m(\om_0)\le (\Phi(\om_0))(s_0) + \tfrac13\e ,$$ and
$$ (\Phi(\om_0))(s_0) < (\Phi(\om))(s_0) +\tfrac13\e$$
for all $\om\in O_0$. Hence, for $\mu$-almost all $\om\in K\cap O_0$,
$$
m(\om) -\tfrac13\e \le M - \tfrac13\e \le
m(\om_0)\le (\Phi(\om_0))(s_0) + \tfrac13\e < (\Phi(\om))(s_0) + \tfrac23\e.
$$
It follows that the Borel set $(K\cap O_0)\setminus N$, where $N$ is some Borel set satisfying $\mu(N)=0$, may be added to the collection $(\Om_i)_{i\in I}$. This contradicts the maximality of this family.
\end{proof}

\begin{corollary}\label{cor3}
Let $\mu$ be a finite Radon measure on a topological space $\Om$
and let $X$ and $Y$ be Banach spaces. Consider a function $\Phi: \Om\to \calL(X,Y)$.
\ben
\item If $ \om\mapsto \n\Phi(\om)x\n$ is lower
semicontinuous for all $x\in X$, then $\Phi$ is strongly $\mu$-normable.
\item If $ \om\mapsto | \langle \Phi(\om)x, y^* \rangle|$ is lower
semicontinuous for all $x\in X$ and $y^*\in Y^*$, then $\Phi$ is weakly $\mu$-normable.
\een
\end{corollary}

Here are two further examples.

\begin{example} \label{ex3}
Consider $\Omega=(0,1)$, $X=L^1(0,1)$, $Y=\K$,  and let $\Phi:(0,1)\to
L^\infty(0,1) = \calL(L^1(0,1),\K)$ be defined by $\Phi(t) :=
1_{(0,t)}$. For all $x\in L^1(0,1)$ the function $t\mapsto
\Phi(t)x = \int_0^t x(s)\,ds$ is continuous. Corollary \ref{cor3}
asserts that $\Phi$ is normable. In fact, for $f(t) :=
\frac1t 1_{(0,t)}$ one even has
$$ \n \Phi(t)\n = |\Phi(t)f(t)| = 1 \quad \forall t\in (0,1).$$
\end{example}

\begin{example} \label{ex4} Let $X_1,X_2$ be Banach spaces and let $T:X_1\to
X_2$ be a bounded linear operator with  $\|T\|=1$. Consider
$\Omega=[0,1]$, $X= C([0,1],X_1)$, $Y=X_2$ and let $\Phi:\Om\to\L(X,Y)$ be
defined by $\Phi(t):=T_t$, where $T_t(x)= T(x(t))$ for $x\in X$. For
all $x\in X$ the function $t\mapsto T_t x$ is continuous.
Corollary \ref{cor3} asserts  that $\Phi$ is weakly (and hence strongly) normable. In
fact, for each $\e>0$ and $t\in [0,1]$ we can select $x^\e\in
B_{X_1}$ and $y^{*\e}\in B_{X_2\s}$ such that $|\lb T x^\e, y^{*\e}\rb|> 1-\e$. Defining
$f^\e := 1\otimes x^\e$ and $g^\e := 1\otimes y^{*\e}$ one  has
$$ \n \Phi(t)\n \le |\lb \Phi(t)f^\e(t), g^\e(t)\rb| +\e\quad \forall t\in [0,1].$$
\end{example}

In the Examples  \ref{ex1}, \ref{ex2} and \ref{ex3}
the norming was exact. The next proposition formulates a simple
sufficient (but by no means necessary) condition for this to be possible:

\begin{proposition}
Let $X$ and $Y$ be Banach spaces and consider a function $\Phi:\Om\to\L(X,Y)$.
\ben
\item Suppose that $\Phi:\Om\to \calL(X,Y)$ is strongly $\mu$-normable. If $X$ is reflexive, there exists a strongly $\mu$-measurable function $f:\Om\to X$ such that
for $\mu$-almost all $\om\in \Om$ one has $\n f(\om)\n \le 1$ and
$$ \n \Phi(\om)\n = \n \Phi(\om)f(\om)\n.$$
\item Suppose that $\Phi:\Om\to \calL(X,Y)$ is weakly $\mu$-normable. If $X$ and $Y$ are reflexive, there exist strongly $\mu$-measurable functions $f:\Om\to X$
and $g:\Om\to Y\s$ such that
for $\mu$-almost all $\om\in \Om$ one has $\n f(\om)\n \le 1$, $\n g(\om)\n\le 1$, and
$$ \n \Phi(\om)\n = |\lb\Phi(\om)f(\om), g(\om)\rb|.$$
\een
\end{proposition}
\begin{proof}
We shall prove (1), the proof of (2) being similar.

For every $n\ge 1$ choose a strongly $\mu$-measurable function $f_n:\Om\to X$ such that
for $\mu$-almost all $\om\in \Om$ one has $\n f_n(\om)\n \le 1$ and
$$ \n \Phi(\om)\n \le \n \Phi(\om)f_n(\om)\n +\tfrac1n.$$
Since $\mu$ is finite, the sequence $(f_n)_{n=1}^\infty$ is bounded in the reflexive space $L^2(\mu;X)$ and therefore it has a weakly convergent subsequence $(f_{n_k})_{k=1}^\infty$. Let $f$ be its weak limit. By Mazur's theorem there exist convex combinations $g_j$ in the linear span of $(f_{n_k})_{k=j}^\infty$ such that $\n g_j-f\n < \frac1j$.
By passing to a subsequence we may assume that $\lim_{j\to\infty} g_j = f$ $\mu$-almost surely. Clearly, for $\mu$-almost all $\om\in\Om$ one has
$\n g_j(\om)\n\le 1$ and
$$ \n \Phi(\om)\n \le \n \Phi(\om)g_j(\om)\n +\tfrac1{n_j}.$$
The result follows from this by passing to the limit $j\to\infty$.
\end{proof}

The following example shows that the separability condition of Proposition \ref{prop} and the lower semicontinuity assumption of Proposition \ref{prop:BA} and its corollaries cannot be omitted, even when $X$ is a Hilbert space and $Y=\K$.

\begin{example}\label{ex}
Let $\Om=(0,1)$, $X = l^2(0,1)$, and $Y=\K$. Recall that
$l^2(0,1)$ is the Banach space of all functions $\phi:(0,1)\to \R$
such that $$\n \phi\n^2 := \sup_{U\in \mathscr{U}} \big\{\sum_{t\in U} |\phi(t)|^2\big\} <\infty,$$
where $\mathscr{U}$ denotes the set of all finite subsets of $(0,1)$.
Note that for all $\phi\in l^2(0,1)$ the set of all $t\in (0,1)$ for which $\phi(t)\not=0$ is at most countable; this set will be referred to as the {\em support} of $\phi$.

Define $\Phi:(0,1)\to \calL(l^2(0,1),\K)$ by
$$ \Phi(t)\phi:=
\phi(t).$$
Clearly, $\n \Phi(t)\n=1$ for all $t\in (0,1)$.
Also, $\Phi(t)\phi = 0$ for all $t$ outside the countable support of $\phi$ and therefore this function is always measurable.

Suppose now that a strongly measurable function $f:(0,1)\to l^2(0,1)$
exists such that $$ 1 \le |\Phi(t)f(t)|+\tfrac12$$ for almost all $t\in (0,1)$.
Let $N$ be a null set such that this inequality holds for all $t\in (0,1)\setminus N$.
For $t\in (0,1)\setminus N$ it follows that $|(f(t))(t)| \ge \frac12$.
Let $f_n:(0,1)\to l^2(0,1)$ be simple functions such that $\limn f_n = f$ pointwise almost everywhere, say on $(0,1)\setminus N'$ for some null set $N'$. The range of each $f_n$ consists of finitely many elements of $l^2(0,1)$, each of which
has countable support. Therefore there exists a countable set $B\subseteq(0,1)$ such that the support of $f(t)$ is contained in $B$
for all $t\in (0,1)\setminus N'$.
For $t\in (0,1)\setminus (N\cup N')$, the inequality
$|(f(t))(t)| \ge \frac12$ implies that $t\in B$. Hence, $(0,1)\setminus(N\cup N')\subseteq B$, a contradiction.
\end{example}

\section{Spaces of operator-valued functions}

Throughout this section, $(\Om,\Sigma,\mu)$ is a finite measure space and $X$ and $Y$ are Banach spaces.

We introduce the linear spaces
$$
\bal {\M}(\mu;\L(X,Y))  &:= \{\Phi:\Om\!\to\! \L(X,Y):
 \text{$\Phi $ is strongly $\mu$-measurable }\},\\
{\M}\sot(\mu;\L(X,Y))  &:= \{\Phi:\Om\!\to\! \L(X,Y):
 \text{$\Phi x$ is strongly $\mu$-measurable $\forall x\in X$}\},\\
{\M}\wot(\mu;\L(X,Y))  &:= \{\Phi:\Om\!\to\! \L(X,Y):  \text{$\Phi x$ is
weakly $\mu$-measurable $\forall x\in X$}\}.
\eal
$$
Two functions $\Phi_1$ and $\Phi_2$ in ${\M}(\mu;\L(X,Y))$ are
identified when $\Phi_1 = \Phi_2$ $\mu$-almost everywhere, two functions
$\Phi_1$ and $\Phi_2$ in ${\M}\sot (\mu;\L(X,Y))$ are identified
when $\Phi_1x = \Phi_2x$ $\mu$-almost everywhere for all $x\in X$,
and  $\Phi_1$ and $\Phi_2$ in ${\M}\wot(\mu;\L(X,Y))$
are identified when $\lb \Phi_1x,y\s\rb = \lb \Phi_2x,y\s\rb$
$\mu$-almost everywhere for all $x\in X$ and $y\s\in Y\s$.

As special cases, for $X=\K$  we put $ \M(\mu;X)  :=
{\M}(\mu;\L(\K,X))$ (which coincides with $
{\M}\sot(\mu;\L(\K,X))$) and $\M\wot(\mu;X) :=
{\M}\wot(\mu;\L(\K,X))$.

The following easy fact will be useful below.

 \begin{proposition} \label{meas} For $\Phi\in {\M}\sot(\mu;\L(X,Y))$ and $f\in \M(\mu;X)$, $$g(\om):=\Phi(\om)f(\om)$$ defines a function $g\in \M(\mu;Y)$.
 \end{proposition}
 \begin{proof}
For simple functions $f$ this is clear. The general case follows from this, using that $\mu$-almost everywhere limits of strongly $\mu$-measurable functions are
strongly $\mu$-measurable.
 \end{proof}

For $1\le p\le \infty$ we consider the normed linear spaces
$$
\bal L^p(\mu;\L(X,Y)) & :=\Big\{ \Phi\in \M(\mu;\L(X,Y)): \
\|\Phi\|_{L^p(\mu;\L(X,Y))} <\infty\Big\},\\
L^p\sot(\mu;\L(X,Y)) & :=\Big\{ \Phi\in \M\sot(\mu;\L(X,Y)): \
\|\Phi\|_{L\sot^p(\mu;\L(X,Y))} <\infty\Big\},\\
L\wot^p(\mu;\L(X,Y)) & :=\Big\{ \Phi\in \M\wot(\mu;\L(X,Y)): \
\|\Phi\|_{L\wot^p(\mu;\L(X,Y))}<\infty\Big\},
\eal
$$
where
$$
\bal \|\Phi\|_{L^p(\mu;\L(X,Y))} & := \Big(\int_\Omega \|
\Phi(\om)\|^p\, d\mu(\om)\Big)^{1/p}, \\
\|\Phi\|_{L\sot^p(\mu;\L(X,Y))} & := \sup_{\n x\n\le 1}
\Big(\int_\Omega \|
\Phi(\om)x\|^p\, d\mu(\om)\Big)^{1/p}, \\
\|\Phi\|_{L\wot^p(\mu;\L(X,Y))} & := \sup_{\n x\n\le 1} \sup_{\n y\s\n\le 1}\Big(\int_\Omega |\lb \Phi(\om)x,y\s\rb|^p \, d\mu(\om)\Big)^{1/p},\eal
$$
with the obvious modifications for $p=\infty$. As special cases we
write
$L^p(\mu;X):= L^p(\mu;\L(\K,X))=L\sot^p(\mu;\L(\K,X))$ and
$L^p\wot(\mu;X):= L^p\wot(\mu;\L(\K,X))$.
Note that all these definitions  agree with the usual ones.

Let us recall some spaces of vector measures that are used in the
sequel. The reader is referred to  \cite{DieUhl} and \cite{Din} for the concepts needed in this paper. Fix $1\le p\le \infty$ and let $E$ be a Banach space. We denote by $V^p(\mu;E)$ the
Banach space of all vector measures $F:\Sigma\to E$ for which
$$\|F\|_{V^p(\mu;E)}:= \sup_{\pi\in \parti(\Om)} \Big\n\sum_{A\in \pi} \frac{1}{\mu(A)}(1_A\otimes F(A)) \Big\n_{L^p(\mu;E)}
<\infty,
$$
where $\parti (\Om)$ stands for the collection of all finite
partitions of $\Om$ into disjoint sets of strictly positive $\mu$-measure. Similarly we denote by $V\wot^p(\mu;E)$ the Banach
spaces of all vector measures $F:\Sigma\to E$ for which
$$\|F\|_{V\wot^p(\mu;E)}:= \sup_{\pi\in \parti(\Om)}\Big\n\sum_{A\in \pi} \frac{1}{\mu(A)}(1_A\otimes F(A)) \Big\n_{L\wot^p(\mu;E)}
<\infty.
$$
In both definitions of the norm we make the obvious modification
for $p=\infty$. Note that $\|F\|_{V^1(\mu;E)}$ and
$\|F\|_{V\wot^1(\mu;E)}$ equal the variation and semivariation of
$F$ with respect to $\mu$, respectively. It is well known that for
$1\le p<\infty$ and $1/p+1/{p'}=1$ one has a natural isometric isomorphism
$$(L^p(\mu;E))^*\simeq V^{p'}(\mu; E^*).$$

We now concentrate on the case $E=\L(X,Y)$. For each $\Phi\in
L^{1}( \mu;\L(X,Y))$ one may define a vector measure $F:\Sigma\to \L(X,Y)$ by $$F(A):=\int_A
\Phi\,d\mu$$ which satisfies
$$\|F\|_{V^1(\mu;\L(X,Y))}= \|\Phi\|_{L^{1}( \mu;\L(X,Y))}.$$
In the next proposition we extend this definition to functions
$\Phi\in L\sot^{p}( \mu;\L(X,Y))$, $1<p<\infty$. The case $p=1$ will be addressed in Remark \ref{rem} and Theorem \ref{t1}.

\begin{proposition} \label{vp} Assume that $\Phi\in L\sot^p( \mu;\L(X,Y))$ for some $1<p<\infty$. Define $F:\Sigma\to\L(X,Y)$ by
$$F(A)x := \int_A \Phi(\om)x\,d\mu(\om),\quad x\in X.$$
Then $F$ is an $\L(X,Y)$-valued vector measure and, for any $q\in [1,p]$, one has
$$\n F\n_{V\wot^q(\mu;\L(X,Y))}
\le \|\Phi\|_{L\sot^q( \mu;\L(X,Y))}.$$
\end{proposition}
\begin{proof}
Let us first prove that $F$ is countably additive. Let $(A_n)_{n\ge 1}$
be a sequence of pairwise disjoint sets in $\Sigma$ and let
$A=\bigcup_{n\ge 1}A_n$. Put $T:=F(A)$ and $T_n:=F(A_n)$. Then,
$$\bal
\Big\|T-\sum_{n=1}^N T_n\Big\|&= \sup_{\|x\|=1}\Big\|Tx-\sum_{n=1}^N T_nx\Big\|\\
&=\sup_{\|x\|=1} \Big\|\int_{\bigcup_{n\ge N+1} A_n} \Phi(\om) x \,d\mu(\om)\Big\|\\
&\le \sup_{\|x\|=1} \Big(\int_{\Omega} \|\Phi(\om) x\|^p\,d\mu(\om)\Big)^{1/p}\mu\Big(\bigcup_{n\ge N+1} A_n\Big)^{1/p'}\\
&\le \|\Phi\|_{L^p( \mu;\L(X,Y))}\ \mu\Big(\bigcup_{n\ge N+1} A_n\Big)^{1/p'}.
\eal
$$
Hence $T=\sum_{n\ge 1} T_n$ in $\L(X,Y)$. Next,
$$\bal
\n F\n_{V\wot^q(\mu;\L(X,Y))} & = \sup_{\pi\in \parti (\Omega)}
\sup_{\|e^*\|=1}  \Big(\sum_{A\in \pi}
\frac{|\langle F(A),e^*\rangle|^q}{(\mu(A))^{q-1}}\Big)^{1/q}\\
& = \sup_{\pi\in \parti (\Omega)} \sup_{\|e^*\|=1}
\sup_{\|(\a_A)\|_{q'}=1}
\Big|\sum_{A\in \pi} \alpha_A \big\langle\frac{F(A)}{(\mu(A))^{1/q'}}, e^*\big\rangle\Big|\\
& =  \sup_{\pi\in \parti(\Omega)}\sup_{\|(\a_A)\|_{q'}=1}
\Big\|\sum_{A\in \pi}
\alpha_A \frac{F(A)}{(\mu(A))^{1/q'}}\Big\|_{\L(X,Y)} \\
& =  \sup_{\pi\in
\parti(\Omega)}\sup_{\|(\a_A)\|_{q'}=1}\sup_{\|x\|=1}
\Big\|\sum_{A\in
\pi} \alpha_A \frac{F(A)}{(\mu(A))^{1/q'}}x\Big\| \\
& = \sup_{\pi\in \parti(\Omega)}\sup_{\|(\a_A)\|_{q'}=1}\sup_{\|x\|=1} \Big\|\int_\Omega \Big(\sum_{A\in \pi} \alpha_A \frac{1_A}{(\mu(A))^{1/q'}}\Big) \Phi(\om)x \,d\mu(\om)\Big\|\\
& \le  \sup_{\|x\|=1} \Big(\int_\Omega \| \Phi(\om)x \|^q\,d\mu(\om)\Big)^{1/q}\\
& =\|\Phi\|_{L\sot^{q}(\mu;\L(X,Y))}. \eal
$$
\end{proof}

\begin{remark}\label{rem} The same results holds for functions $\Phi\in L\sot^1(
\mu;\L(X,Y))$ provided the family $\{\om\mapsto\Phi(\om)x: x\in B_X\}$ is equi-integrable in $L^1(\mu;X)$.
\end{remark}

The next definition introduces a new class of Banach spaces intermediate  between $L^p(\mu;\L(X,Y))$ and $L^p\sot(\mu;\L(X,Y))$.

\begin{definition} For $1\le p\le \infty$ we consider the Banach space
$$L^{p}[\mu;\L(X,Y)]:=\{\Phi\in
\M\sot(\mu;\L(X,Y)): \ \n \Phi\n_{L^{p}[
\mu;\L(X,Y)]}<\infty\},$$ where
$$
 \n \Phi\n_{L^{p}[
\mu;\L(X,Y)]} := \sup_{\|f\|_{L^{p'}(\mu;X)}=1}\int_\Omega \|
\Phi(\om)f(\om)\|\,d\mu(\om).$$
\end{definition}
It is clear that
$$ L^p(
\mu;\L(X,Y))\embed L^p[ \mu;\L(X,Y)]\embed
L\sot^p(\mu;\L(X,Y))$$ with contractive inclusion mappings.
Using these spaces we can  prove the following improvement of Proposition
\ref{vp}.
\begin{theorem} \label{t2} Let $1< p<\infty$. Then
$$L^p[\mu;\L(X,Y)]\embed V^p(\mu; \L(X,Y))$$
and the inclusion mapping is contractive.
\end{theorem}
\begin{proof}   Using the inclusion
into $L^p[\mu;\L(X,Y)]\embed L\sot^p( \mu;\L(X,Y))$, from Proposition
\ref{vp} we see that $F(A)x:=\int_A \Phi(\om)x\,d\mu(\om)$ defines a vector measure $F:\Sigma\to\L(X,Y)$.

Now, if $\pi\in\parti(\Om)$, then for $\e>0$ and each $A\in \pi$
there exist $x_A\in B_X$ and $y_A^*\in B_{Y^*}$ so that
$$\|F(A)\|^p< \Big|\Big\langle \int_A \Phi(\om)x_A \,d\mu(\om), \,y_A^*\Big\rangle\Big|^p+ \frac{\e}{\hbox{card}(\pi)}.$$
Hence,
$$\bal
\ & \sum_{A\in\pi} \frac{\|F(A)\|^p}{(\mu(A))^{p-1}}
\\ & \quad \le \sum_{A\in\pi}
\frac1{(\mu(A))^{p-1}}\Big|\Big\langle \int_A \Phi(\om)x_A\, d\mu(\om), \,y_A^*\Big\rangle\Big|^p + \e \\
& \quad \le \sup_{\|(\beta_A)\|_{p'}=1}\Big(\sum_{A\in\pi}
\frac1{(\mu(A))^{1/p'}}\Big\langle \int_A \Phi(\om)x_A\, d\mu(\om), \,\beta_Ay_A^*\Big\rangle\Big)^p + \e \\
& \quad \le \sup_{\|(\beta_A)\|_{p'}=1}\Big(\int_\Omega
\Big\langle  \Phi(\om)\sum_{A\in\pi} 1_A\otimes\frac{\beta_A x_A}{(\mu(A))^{1/p'}},  \sum_{A\in\pi}1_A \otimes y_A^* \Big\rangle\,d\mu(\om)\Big)^p +\e \\
& \quad \le \sup_{\|(\beta_A)\|_{p'}=1}\Big(\int_\Omega \Big\|
\Phi(\om)\sum_{A\in\pi} 1_A\otimes \frac{\beta_A x_A}{(\mu(A))^{1/p'}}\Big\|
\,d\mu(\om)\Big)^p + \e \\
& \quad\le \sup_{\|f\|_{L^{p'}(\mu;X)}=1}\Big(\int_\Omega
\| \Phi(\om)f(\om)\| \,d\mu(\om)\Big)^p + \e \\
& \quad \le \|\Phi\|^p_{L^p[\mu;\L(X,Y)]}+\e. \eal
$$
Since $\e>0$ was arbitrary, this gives the result.
 \end{proof}

For $1\le p,q<\infty$ we define
$${\rm Mult\/}[L^p(\mu;X),L^q(\mu;Y)]$$ to be the linear space of all $\Phi\in \M\sot(\mu;\L(X,Y))$ such that $\om\mapsto\Phi(\om)f(\om)$ belongs to
$L^q(\mu;Y)$ for all $f\in L^p(\mu;X).$
By a closed graph argument the linear operator $M_\Phi: f\mapsto
\Phi f$ is bounded, and the space ${\rm
Mult\/}[L^p(\mu;X),L^q(\mu;Y)]$ is a Banach space with respect to
the norm $$\|\Phi\|_{{\rm
Mult\/}[L^p(\mu;X),L^q(\mu;Y)]}:=\|M_\Phi\|_{\L(L^p(\mu;X),L^q(\mu;Y))}.$$
We refer to \cite{FouSch} for further details and and some results on 
spaces of  multipliers between different spaces of vector valued functions, extending those proved in \cite{ArreBlas} for sequence spaces.

\begin{theorem}\label{thm:mult} Let $X$ and $Y$ be Banach spaces and let $1\le q<p<\infty$. We have a natural isometric isomorphism
$$ {\rm Mult\/}[L^{p}(\mu;X),L^q(\mu;Y)]\simeq L^{r}[ \mu;\L(X,Y)],$$
where $\frac1r=\frac1q-\frac1p$.
\end{theorem}
\begin{proof} The case $q=1$ corresponds to $r=p'$ and the result is just the definition of the space $L^{p'}[ \mu;\L(X,Y)]$.
Assume $1<q<p$ and $\Phi\in L^{r}[ \mu;\L(X,Y)]$.

Let $f\in L^{p}(\mu;X)$.  Then for any $\phi\in {L^{q'}(\mu)}$ we
have that $\om \to  f(\om)\phi(\om)$ belongs to $L^{r'}(\mu;X)$.
Hence
$$\int_\Om \|\Phi(\om)f(\om)\| |\phi(\om)|\,d\mu(\om)\le
\|\Phi\|_{ L^{r}[ \mu;\L(X,Y)]}\|\phi\|_{L^{q'}(\mu)}
\|f\|_{L^{p}(\mu;X)}.$$ Taking the supremum over the unit ball of
${L^{q'}(\mu)}$ the first inclusion is achieved.

Conversely, let $\Phi\in {\rm Mult\/}[L^{p}(\mu;X),L^q(\mu;Y)]$.
Let $g\in L^{r'}(\mu;X)$, and choose $\psi\in
L^{q'}(\mu)$ and $f\in L^{p}(\mu;X)$ in such a way that $g=\psi f$ and
$$\|g\|_{L^{r'}(\mu;X)}=\|\psi\|_{L^{q'}(\mu)}\|f\|_{L^{p}(\mu;X)}.$$
Now observe that $\Phi(\om)g(\om) = \psi(\om) \Phi(\om)f(\om)\in
L^1(\mu;Y)$ and
$$\int_\Omega\|\Phi(\om)g(\om)\|\,d\mu(\om)\le
\|\psi\|_{L^{q'}(\mu)}\|\Phi\|_{{\rm
Mult\/}[L^{p}(\mu;X),L^q(\mu;Y)]}\|f\|_{L^{p}(\mu;X)}.$$
\end{proof}

The next result establishes a link with the notion of strong $\mu$-measurability.

 \begin{proposition}  \label{inclusion} Let $X$ be a Banach space,  let $1\le p\le \infty$, and let $\Phi\in
L^p[ \mu;\L(X,Y)]$ be strongly $\mu$-normable. Then $\om\mapsto\n \Phi(\om)\n$ belongs to $ L^p(\mu)$ and
$$\Big(\int_\Omega \|\Phi(\om)\|^p\,d\mu(\om)\Big)^{1/p}\le \|\Phi\|_{L^p[
\mu;\L(X,Y)]}.$$
\end{proposition}
\begin{proof} By assumption, for any $\e>0$ there
exists $f^\e\in \M(\mu;X)$ such that for $\mu$-almost all $\om\in
\Omega$ one has $\|f^\e(\om)\|\le 1$ and $\|\Phi(\om)\|\le
\|\Phi(\om)(f^\e(\om))\|+\e.$ 

If $\e_n\downarrow 0$, then for $\mu$-almost all $\omega\in \Omega$  
$$\n \Phi(\omega)\n =  
\lim_{n\to\infty} \n \Phi(\omega)f^{\e_n}(\omega)\n.$$
The strong $\mu$-measurab\-ility of $\omega\mapsto \Phi(\omega)x$ for all $x\in X$ implies the 
 the strong $\mu$-measurability of the functions $\omega\mapsto \Phi(\omega)f^{\e_n}(\omega)$. 
It follows that $\omega\mapsto \n \Phi(\omega)\n$ is $\mu$-measurable. 

Let $\phi\in L^{p'}(\mu)$ and consider
$\om\to \phi(\om)f^\e(\om)\in L^{p'}(\mu;X)$. Then
$$
\bal \int_{\Omega} \|\Phi(\om)\||\phi(\om)|\,d\mu(\om) &\le
\int_\Omega
\|\Phi(\om)(\phi(\om)f^\e(\om))\|\,d\mu(\om) +\e\|\phi\|_{L^1(\mu)}\\
&\le \|\Phi\|_{L^p[ \mu;\L(X,Y)]}\|\phi\|_{L^{p'}(\mu)}
+\e\|\phi\|_{L^1(\mu)}. \eal
$$
Since $\e>0$ was arbitrary, this gives the result.
 \end{proof}

By invoking Proposition \ref{prop} we shall now deduce some further results
under the assumption that the space $X$ is separable.
The first should be compared the remarks preceding Proposition \ref{vp} (where functions
$\Phi\in L^1(\mu;\L(X,Y))$ are considered) and Remark \ref{rem} (where functions
$\Phi\in L\sot^1(\mu;\L(X,Y))$ are considered).
\begin{theorem} \label{t1}
 Let $X$ be a separable Banach space and let $\Phi\in L^{1}[ \mu;\L(X,Y)]$ be given. Define $F:\Sigma\to\L(X,Y)$ by
$$F(A)x := \int_A \Phi(\om) x\,d\mu(\om), \quad x\in X.$$
Then $F$ is an $\L(X,Y)$-valued vector measure and
$$\n F\n_{V^1(\mu;\L(X,Y))}
\le \|\Phi\|_{L^1[ \mu;\L(X,Y)]}.$$
\end{theorem}
\begin{proof} First we prove that $F$ is countably additive. Let $(A_n)_{n\ge 1}$ be a sequence of pairwise disjoint sets in $\Sigma$ and let $A=\bigcup_{n\ge 1} A_n$. Put $T:=F(A)$ and $T_n:=F(A_n)$.
Combining Proposition \ref{prop} and Proposition \ref{inclusion} one obtains that $\|\Phi\|\in L^1(\mu)$. Hence,
$$\bal
\Big\|T-\sum_{n=1}^N T_n\Big\|&= \sup_{\|x\|=1}\Big\|Tx-\sum_{n=1}^N T_n x\Big\|\\
&=\sup_{\|x\|=1} \Big\|\int_{\bigcup_{n\ge N+1} A_n} \Phi(\om) x \,d\mu(\om)\Big\|\\
&\le \int_{\bigcup_{n\ge N+1} A_n} \|\Phi(\om)\|\,d\mu(\om).
\eal
$$
Hence $T=\sum_{n\ge 1} T_n$ in $\L(X,Y)$. Next, using that $\|F(A)\|\le \int_A \|\Phi(\om)\|\,d\mu(\om)$, from Proposition \ref{inclusion} we conclude that
$$
\n F\n_{V^1(\mu;\L(X,Y))}
 = \sup_{\pi\in\parti(\Om)} \sum_{A\in \pi}\n F(A)\n
 \le \|\Phi\|_{L^{1}[\mu;\L(X,Y)]}.
$$
\end{proof}

Our final result extends the factorization result that was used in the proof of Theorem \ref{thm:mult}.

\begin{theorem}\label{fact} Let $1\le p_1,p_2, p_3 <\infty$ satisfy $\frac1{p_1}=\frac1{p_2}+\frac1{p_3}$ and let $X$ be a separable Banach space. A function
$\Phi\in \M\sot(\mu;\L(X,Y))$ belongs to $L^{p_1}[\mu;\L(X,Y)]$
if and only if $ \Phi = \psi\Psi$ for suitable functions $\psi\in
L^{p_2}(\mu)$ and $\Psi\in L^{p_3}[\mu;\L(X,Y)]$. In this
situation we may choose $\psi$ and $\Psi$ in such a way that
$$ \n \Phi\n_{L^{p_1}[\mu;\L(X,Y)]} = \n \psi\n_{L^{p_2}(\mu)} \n \Psi\n_{L^{p_3}[\mu;\L(X,Y)]}.$$
\end{theorem}

\begin{proof}  
To prove the `if' part let  $\Phi\in L^{p_1}[\mu;\L(X,Y)]$. Using Proposition \ref{inclusion} together with Proposition \ref{prop}
one has that $\|\Phi\|\in L^{p_1}(\mu)$. Put
$$\psi(t) := \|\Phi(t)\|^{p_1/p_2}, \quad
\Psi(t)  := \left\{\begin{array}{ll}
\|\Phi(t)\|^{p_1/p_3}
\frac{\Phi(t)}{\|\Phi(t)\|}& \text{if}\ \Phi(t)\ne 0, \\
0 & \text{if}\ \Phi(t)=0.
\end{array}
\right.
$$
Clearly $\psi\in L^{p_2}(\mu)$ and $\Psi\in L^{p_3}[\mu;\L(X,Y)]$. Now 
 for each $g\in L^{p_3'}(\mu;X)$, invoking Proposition \ref{meas}, one has that $\Psi g \in \M(\mu, Y)$ and $$\|\Psi(t)g(t)\|
\le
\|\Phi(t)\|^{p_1/p_3}\|g(t)\|.$$
Hence the right hand side defines a function in $L^1(\mu)$ and therefore $\Psi g\in L^1(\mu,Y)$.
The above decomposition satisfies the required identity for the norms.

To prove the `only if' part let $\psi\in
L^{p_2}(\mu)$ and $\Psi\in L^{p_3}[\mu;\L(X,Y)]$ be given. For each  $f\in L^{p'_1}(\mu;X)$  we have $\psi f\in
L^{p'_3}(\mu; X)$. Hence $\Psi(\psi f)\in L^1(\mu;Y)$.

\end{proof}


\end{document}